\documentclass[12pt]{article}
\usepackage{graphicx}
\usepackage{amssymb}
\usepackage{amsmath}
\usepackage{amsthm}
\usepackage{mathscinet}
\usepackage{a4}
\usepackage{enumerate}
\usepackage{url}

\usepackage[usenames, dvipsnames]{color}
\newcommand{\corr}[1]{{\textcolor{black}{#1}}}

\newcommand{\HA}{\ensuremath{H^2(\mathbb{A})}}

\theoremstyle{plain}

\usepackage{thmtools}

\declaretheoremstyle[%
  spaceabove={\topsep},%
  spacebelow={\topsep},%
  bodyfont={\itshape}, %
  headfont={\bfseries},%
  postheadspace=1em,%
  qed=\qedsymbol,%
  headpunct={}
]{mystyle2}

\theoremstyle{mystyle2}
\newtheorem*{theorem*}{} 

\definecolor{cb-blue-green}{RGB}{0,073,073}
\definecolor{cb-blue}{RGB}{0,109,219}

\newcommand{\corrn}[1]{#1}

\begin{document}

\title{Dynamic mode decomposition for analytic maps}

\author{Julia Slipantschuk, Oscar F. Bandtlow, and Wolfram Just\\
School of Mathematical Sciences,\\ Queen Mary University of London,
London, UK}

\date{May 13, 2019}

\maketitle

\begin{abstract}
Extended dynamic mode decomposition (EDMD) provides a class
of algorithms to identify patterns and effective degrees of
freedom in complex dynamical systems.
We \corr{show}
that the modes identified by
EDMD correspond to those of compact \corrn{Perron-Frobenius and Koopman}
operators defined on suitable Hardy-Hilbert spaces
when the method is applied to classes of analytic maps.
Our findings elucidate the interpretation
of the spectra obtained by EDMD for complex dynamical systems.
We illustrate our results by numerical simulations
for analytic maps.
\end{abstract}

\section{Introduction}

The quest to identify effective degrees of freedom in a complex dynamical
system is a fundamental topic in almost all branches of science.
The archetype and historical origin of this endeavour can be seen in
the derivation of thermodynamics from microscopic equations of
motion within a hydrodynamic description. Here, the relevant
macroscopic densities are determined by the classical conservation
laws of physics. In a mathematical setting the problem of identifying
effective degrees of freedom and reducing the dynamical description
to a lower dimensional set of equations can be cast in terms of
centre manifold reductions \cite{GuHo:86} or adiabatic elimination
procedures \cite{vKam_PR85}. In these classical situations
the reduction of the description to effective degrees of freedom has resulted in
the derivation of transport equations for systems far from equilibrium
using projection operator techniques \cite{Naka_PTP58,Zwan_JCP60},
an understanding
of how dissipation emerges in many particle
Hamiltonian systems \cite{Prigo}, or the Bandtlow-Coveney equation
for transport properties in discrete-time dynamical systems \cite{BCov}, to name but a few.
With the advent of the study of complex and chaotic dynamical behaviour the focus
has shifted and broadened. Nowadays the problem of identifying relevant
degrees of freedom occupies a diverse range of scientific fields
ranging from physics
and the life sciences, to computational studies of pattern recognition or data processing. Many of the current
algorithmic approaches have been inspired by the classical ideas mentioned above.
In a nutshell, these methods are based on identifying an optimal mode decomposition which can be
used to effectively describe the system of interest.

A recent instalment of these ideas has become known
as Dynamic Mode Decomposition, introduced in
\cite{Sch10} and extended in \cite{WKR15}; see also \cite{KlKoSc_JCD16}
for an illustration of this concept in the context of nonlinear
stochastic systems.
At their core, these methods condense the dynamical observations into a suitably chosen effective
linear evolution matrix.
The eigenvalues and eigenvectors
of this matrix then provide information concerning the relevance and structure of
the effective degrees of freedom of the system, and can sometimes be related to
spectral data of global evolution operators, known as Perron-Frobenius operators, or their formal adjoints, known as Koopman operators.
Data-driven methods to approximate Perron-Frobenius and Koopman operators have
become prominent with the work of  Dellnitz and Junge \cite{DelJunge99} and have since been
extended in many ways \cite{DelFroSer00, FroGonQuas14, MeBa04, Gia2017,
KlGePeSc_Nonl18}, to mention but a few.
For an overview of existing algorithms, see the review articles
\cite{Kl_review,BuMoMe_Chaos12, Gia2017} and references therein,
in particular \cite{Kl_review} for comparison of various data-driven algorithms,
\cite{BuMoMe_Chaos12} for an overview of applications in areas of engineering, and the
introductory section of \cite{Gia2017} for a historical overview.

Empirically, these algorithms perform
extremely well and
are essentially fully understood for finite-dimensional linear dynamical
systems. However, open questions remain in more complex dynamical setups, regarding, for
example, under
which conditions the algorithms converge and whether the limiting quantities
are signatures of the underlying dynamical system in the sense that they approximate the spectral
data of the relevant evolution operator.
In this note our aim is to contribute to this issue, by proving
that a certain version of dynamic mode
decomposition, known as Extended Dynamic Mode Decompostion (EDMD), identifies the correct
effective degrees of freedom in an
analytic setup when certain classes of deterministic chaotic dynamical systems are studied.

In order to keep our presentation self-contained we start with
a brief sketch of EDMD in Section~\ref{sec:2}.
In Section~\ref{sec:3} we introduce a class of analytic circle maps, for which rigorous statements
about
EDMD can be made. For this class
of systems we show that EDMD singles out
eigenmodes of the Perron-Frobenius operator 
on a suitably defined space of analytic
functions. In this sense EDMD
effectively performs a coarse graining or smoothing of the dynamics.
We shall also explain that, in an appropriate setting, these results translate into strong spectral
convergence results for the corresponding Koopman operator.
In Section~\ref{sec:4} we illustrate our findings through various
numerical examples based on exactly solvable models. Finally,
in Section~\ref{sec:5}, we put our results in a more general context, including a discussion of
higher dimensional dynamical systems or the relevance of our
rigorous approach for general dynamic mode
decompositions where no proofs can be provided.

\section{Extended Dynamic Mode Decomposition}\label{sec:2}
In the following, we shall provide a brief, informal account of
Extended Dynamic Mode Decomposition (EDMD). Consider a discrete dynamical system
\begin{equation}\label{aa}
z_{n+1}=\tau(z_n) \,
\end{equation}
given by map a $\tau\colon\mathcal{X}\to\mathcal{X}$ on some phase space $\mathcal{X}$.
Assume that the dynamics is observed through a collection
of $N$ scalar functions defined on the phase space given by
$\psi(z) = (\psi_1(z),\ldots,\psi_N(z))^{T}$.
We record the dynamics at a sequence of
$M$ phase space points $z^{(1)},\ldots,z^{(M)}$. These points can be
obtained from a time series if the approach is used as a data
analysis tool, or as a sample from a suitable distribution of points in phase space
if the goal is to investigate the underlying equations of motion \eqref{aa}.
Glossing over details of the underlying theory
(see for example, \cite{WKR15, KlKoSc_JCD16, KlGePeSc_Nonl18, KoMe_JNS18}
and references therein), the fundamental quantity of EDMD is
an $N\times N$ matrix
\begin{equation}\label{eq:A}
A=G H^{-1},
\end{equation} which is constructed from the
observations as follows
\begin{align}
G_{k l} =& \frac{1}{M} \sum_{m=1}^M \psi_k(\tau(z^{(m)}))
\psi_l(z^{(m)})\,, & &(k,l=1,\ldots, N)\,, \label{aba} \\
H_{k l} =& \frac{1}{M} \sum_{m=1}^M \psi_k(z^{(m)})
\psi_l(z^{(m)})\,, & &(k,l=1,\ldots,N) \, . \label{abb}
\end{align}
Given the observations, the matrix $A$ is an optimal representation of the
dynamics in terms of a finite dimensional linear equation of motion in the following sense: it is a
least squares solution to
$AX=Y$ where $X=[\psi(z^{(1)}), \dots, \psi(z^{(M)})]$
and $Y=[\psi(\tau(z^{(1)})), \dots, \psi(\tau(z^{(M)}))]$.\footnote{
A solution to ${\arg\min}_{A}||AX - Y||^2_2$ is given
by $A = Y X^{+}$, where $X^{+}$ is the Moore-Penrose pseudoinverse of $X$.
The matrix $A$ can be written as $A = (Y X^H)(X X^H)^{-1}$ if $XX^H$ is invertible,
where $X^H$ denotes the conjugate transpose of $X$.
Furthermore, $A=(YX^T) (XX^T)^{-1} = G H ^{-1}$, assuming that
for each observable $\psi_i$ there is a $\psi_j$ with $\psi_j(z^{(m)}) =\overline{\psi_i(z^{(m)})}$ for
all $z^{(m)}$, which holds for the observables used in this paper.}
For sufficiently large values of $M$ and $N$
the eigenvalues and eigenvectors of $A$ determine the effective modes
of the system. Eigenmodes with eigenvalues on or close to the complex
unit circle are the slow modes which are relevant for macroscopic behaviour and the long term
dynamics.

The matrix $A$ can be linked to the linear operator governing the
underlying dynamics \eqref{aa}, the Perron-Frobenius operator, or its formal adjoint,
the Koopman operator.
Previous investigations \cite{KoMe_JNS18} have shown
that in the limit of large $N$ and $M$ the matrix $A$ is, in
a certain sense, a suitable matrix representation, provided strong technical conditions are met.
These conditions, however, may be difficult to verify in concrete applications.

Numerical results show
that EDMD can often be applied successfully as a practical algorithm
and indicate that output data (for example, eigenvalues of $A$)
exhibit nice convergence properties.
We will prove that this is indeed the case for a suitable class
of dynamical systems, which will be introduced in the next section.


\section{EDMD for analytic maps}\label{sec:3}
Let us consider a full branch analytic expanding map on an interval, say $[0,2 \pi]$. Using the
canonical mapping
$\varphi\mapsto z=\exp(i\varphi)$ this map can be viewed as a map $\tau$ in the
complex plane leaving the unit circle $\mathbb{T} = \{z\in \mathbb{C}: |z|=1\}$ invariant.
The map $\tau$ will be analytic on $\mathbb{T}$, provided that the branches of the original map on
$[0,2\pi]$ satisfy matching conditions at the endpoints, and will thus have an analytic extension to
an open annulus $\mathbb{A}$ containing $\mathbb{T}$.
This means that such a map admits a Laurent series on
$\mathbb{A}$ which, on the unit circle, coincides with the Fourier
series expansion of the original interval map. Moreover, the Fourier
coefficients will decay exponentially. This last property is one of the
crucial ingredients that will allow us
to define the evolution operators on sufficiently nice function spaces, as discussed below. The
second crucial ingredient is the expansivity of the map, by which we mean that
$|\tau'(z)|>1$ for all $z\in \mathbb{T}$.

The fine statistical properties of the map $\tau$ are captured
by the Perron-Frobenius operator (or transfer operator), which describes
the forward evolution of densities under the action of the system. For analytic expanding
(orientation-preserving) circle maps, it takes the form
\begin{equation}\label{ba}
(\mathcal{L} f)(z)=\sum_j \phi_j'(z) f(\phi_j(z)),
\end{equation}
where $\phi_j$ denotes the $j$-th inverse branch of the analytic map $\tau$.
For instance, for the simple Bernoulli shift map
$\varphi \mapsto 2 \varphi \mbox{ mod } 2 \pi$, the corresponding circle map reads
$\tau(z)=z^2$ with the two inverse branches given by $\phi_1(z)=\sqrt{z}$ and
$\phi_2(z)=-\sqrt{z}$.

The operator in \eqref{ba} is naturally defined on $L^1(\mathbb{T})$,
the positive elements of which are interpretable as probability densities,
but for convenience it is often considered as an operator restricted to $L^2(\mathbb{T})$ so
that Hilbert space methods can be used.
Its adjoint is known as the Koopman operator, which turns out to be
the operator of composition with the map
$\tau$. However,
the Hilbert space $L^2(\mathbb{T})$ is ``fairly large'' so that the spectrum of the bounded operator
$\mathcal{L}$ is the entire closed complex unit disk, with each point in the open unit disk being an
eigenvalue of infinite multiplicity (see, for example, \cite[Remark~4.4]{Kel2}).
Intuitively, the function space $L^2(\mathbb{T})$ simply contains too many ``non-physical''
observables.

In order to capture the
behaviour observed in a time series one often restricts the set of
observables, that is, one takes a suitable subspace of
$L^2(\mathbb{T})$,
so that decay rates show up as isolated spectral points of the Perron-Frobenius operator.
In a sense such a restriction corresponds to the coarse graining used
in statistical mechanics when moving from a conservative microscopic
to a dissipative hydrodynamic description \cite{Prigo}.
An elementary illustration of this aspect
can be found for instance in \cite{SaHa_PLA92a}.\\

\textit{Spectral convergence.}
In order to obtain strong spectral results for our setup of analytic circle
maps, a suitable class of observables is a space of analytic functions.
Following the approach in \cite{BaJuSl_AIHP17}
we restrict observables to be analytic functions on an open annulus $\mathbb{A}$
containing the complex unit circle with an $L^2$-extension to the boundary of $\mathbb{A}$,
the so-called Hardy-Hilbert space \HA.
As shown in \cite{BaJuSl_AIHP17}, the Perron-Frobenius operator given
by \eqref{ba} considered on \HA\ is well-defined and compact, which implies that it
has a discrete spectrum of
eigenvalues which govern the correlation decay and the relaxation of
analytic observables.

In addition, the Perron-Frobenius operator can be effectively approximated by a sequence of
finite rank operators. For this consider an orthogonal basis of the underlying
Hardy-Hilbert space, for instance, the canonical (non-normalised) orthogonal basis
$\psi_m(z)=z^m$ with $m \in \mathbb{Z}$. The corresponding matrix elements
of the operator in \eqref{ba} are then given by
\begin{equation}\label{bb}
L_{k l} = \frac{1}{2 \pi} \int_0^{2 \pi} \psi_k(\tau(\exp(i \varphi)))
\psi_l(\exp(i \varphi))\, d \varphi \, .
\end{equation}

The matrix elements with $k,l=-\bar{N},\ldots,\bar{N}$ yield an
$N \times N$ matrix representation of a finite rank approximation of the Perron-Frobenius operator,
with $N=2 \bar{N}+1$. It has been shown in \cite{BaJuSl_AIHP17} that
these approximations converge to the Perron-Frobenius
operator exponentially fast in operator norm.
Hence, the eigenvalues of the matrices of size $N\times N$
approximate the spectrum of the infinite-dimensional
compact operator ${\mathcal L}$ as $N$ tends to infinity.
Moreover, convergence of the eigenvalues occurs at an exponential rate and
explicit error bounds can be derived from the general theory.
In summary, a finite dimensional matrix approximation using \eqref{bb} provides
the spectrum of the compact Perron-Frobenius operator and that of
its adjoint.

More formally, the results can be stated as follows
\begin{theorem*}
Let $\tau$ be an analytic expanding circe map, and $\mathcal{L}$ the corresponding Perron-
Frobenius operator, given by \eqref{ba}.
Denote by $\{\psi_m\}_{m\in\mathbb{Z}}$ with $\psi_m(z)=z^{m}$ 
the canonical orthogonal basis in \HA and
let $P_N$ be the orthogonal projection operator onto
the subspace spanned by $\psi_{-\bar{N}}, \ldots, \psi_{\bar{N}}$,
where $N = 2\bar{N}+1$.
\begin{enumerate}[(a)]
  \item \textbf{(Compactness of $\mathcal{L}$)}\\
  The operator $\mathcal{L}$ is a well-defined, compact operator from \HA\, to itself.
  \item \textbf{(Matrix representation)} \\
  A matrix representation of the finite rank operator
  $P_N\mathcal{L}P_N$ is given by $(M_{kl})_{k,l}$ with $k,l=-\bar{N},\ldots, \bar{N}$ and
  $M_{kl}=L_{-k, l}$ as in \eqref{bb}.
  \item \textbf{(Convergence in operator norm)}\\
  $||\mathcal{L} - P_N\mathcal{L}P_N||_{\HA\to\HA} = O(\exp(-a N))$ for some $a>0$.
  \item \textbf{(Eigenvalue convergence)}\\
    The spectrum $\sigma(\mathcal{L})$ of $\mathcal{L}$ consists of at most countably many non-zero eigenvalues $\lambda_n$ of finite  multiplicity, with $0$ the only possible accumulation point.
  \begin{enumerate}[(i)]
    \item If
    $(\lambda^{(N)})_{N\in\mathbb{N}}$ with
    $\lambda^{(N)} \in \sigma(P_N\mathcal{L}P_N)$ is a convergent sequence,
    that is, $\lambda^{(N)} \to \lambda$, then $\lambda \in \sigma(\mathcal{L})$.
   \item For every $\lambda \in \sigma(\mathcal{L})$ there exists a sequence
   $(\lambda^{(N)})_{N\in\mathbb{N}}$ with $\lambda^{(N)} \in \sigma(P_N\mathcal{L}P_N)$, such that $\lambda^{(N)} \to \lambda$.
  \end{enumerate}
  More precisely, for suitable enumerations $\lambda_n$ with $n \in \mathbb{N}_0$ of the respective non-zero eigenvalues (taking algebraic multiplicities into account)
  of $P_N \mathcal{L} P_N$ and $\mathcal{L}$, we have: for every $n$
  \begin{equation}\label{eq:eig_conv}
  | \lambda_n(P_N\mathcal{L}P_N) - \lambda_n(\mathcal{L}) |
  = O(\exp(-a N)) \quad \text{as $N\to \infty$}
  \end{equation}
  for some $a>0$.
\end{enumerate}
\end{theorem*}

In practical applications, a-priori error bounds for
\eqref{eq:eig_conv} can be computed explicitly using \cite{Ba08}.

For the proof of the above, observe that \textit{(a)}
follows from \cite[Proposition 3.4]{BaJuSl_AIHP17} and
\textit{(c)} is implied by the proof of \cite[Lemma 3.3]{BaJuSl_AIHP17}. Item \textit{(b)}
follows from a calculation using duality (see, for example, \cite[Lemma 2.3]{SlBaJu_Nonl13}) and
the inner product of \HA\, from
\cite{BaJuSl_AIHP17}.
  Items \textit{(i)} and \textit{(ii)} of \textit{(d)} are known as Properties U and L in \cite{AhLaLi01},
  and follow from Corollaries 2.7 and 2.13 therein, respectively. Finally, the
  exponential convergence of eigenvalues in \eqref{eq:eig_conv} is an immediate
  consequence of \cite[Theorem 2.18]{AhLaLi01} and the ensuing remarks, combined with
  \textit{(c)}.
\\

\textit{The relation between Perron-Frobenius and Koopman operators.}
We have chosen to present our results formulated for the Perron-Frobenius operator.
A large part of the literature on data-driven methods such as EDMD is based on the study of the
Koopman operator, given by $f \mapsto f\circ \tau$, which is the adjoint of the Perron-Frobenius
operator when viewed on $L^2(\mathbb{T})$.
In order to obtain the strong spectral convergence results described above, it was necessary
to restrict the domain of the Perron-Frobenius operator to the ``smaller space" \HA.
This space is densely and continuously embedded in $L^2(\mathbb{T})$, and thus an example of
a test function space 
(see, for example, \cite{SATB95}), 
so that we have

\begin{equation}\label{rigged}
\HA \subset L^2(\mathbb{T}) \simeq L^2(\mathbb{T})' \subset{H^2(\mathbb{A})'},
\end{equation}
where $H^2(\mathbb{A})'$ is the topological dual\footnote{The topological dual $H^2(\mathbb{A})'$
is the space of continuous linear functionals on \HA\, equipped with the topology of uniform
convergence.
Whereas \HA\, consists of analytic functions with exponentially decaying Fourier coefficients, the
space $H^2(\mathbb{A})'$ is ``fairly large'', that is, on top of every function in $L^2(\mathbb{T})$ it
also contains distributions or generalized functions,
with Fourier coefficients allowed to grow exponentially.} of \HA.
The structure \eqref{rigged} is known as a \textit{rigged Hilbert space} or \textit{Gelfand triple} (see, for example, \cite{gelfand_vilenkin} or \cite{bohm_gadella}),
which has been used in the context of dynamical systems to study spectral decompositions for certain chaotic
maps (see, for example, \cite{BAS} and references therein). The (Banach space) adjoint of the Perron-Frobenius
operator $\mathcal{L}$ restricted to
a ``small space" \HA\, can be identified with a Koopman operator extended to a ``large space"
$H^2(\mathbb{A})'$, on which it is compact.

Moreover, it turns out that in our setting of analytic
expanding circle maps, it is even possible to identify
this extended operator on
$H^2(\mathbb{A})'$ with certain Koopman operators 
acting on spaces of analytic functions. As was shown in \cite{BaJuSl_AIHP17}, the space
$H^2(\mathbb{A})'$ can be identified with the space of functions holomorphic on two disks
comprising the complement of the annulus $\mathbb{A}$, denoted $H^2(\mathbb{D}_{in}) \oplus
H_0^2(\mathbb{D}_{out})$. Consequently, the expression \eqref{bb} also yields the matrix
representation of the Koopman operator on this space.\\

\textit{Spectral convergence for EDMD.}
The  results above imply that for analytic circle maps, EDMD
has strong convergence properties and
captures the spectrum of the associated Perron-Frobenius and Koopman operators.
For a suitable choice of sampling points $z^{(m)}$ the expression \eqref{aba}
estimates the matrix elements (\ref{bb}), as established in \cite{WKR15} (see also
\cite{KlKoSc_JCD16} or \cite{KoMe_JNS18}).
For instance, when choosing equidistant points on the unit circle, $z^{(m)}=\exp(2 \pi i m/M)$, the
expression (\ref{aba}) converges exponentially
in $M$ to the integral (\ref{bb}).
The matrix $H$ in \eqref{abb} takes account of the orthonormalisation of
the observables which was used in writing down the matrix elements (\ref{bb}).
Hence, EDMD with equidistant sample points applied to analytic circle
maps with analytic observables
gives precisely the spectrum of the corresponding compact Perron-Frobenius and
Koopman operators. Thus, EDMD singles out the
physically observable decay rates and the corresponding dissipative modes.
We will illustrate this result in the next section by analytically solvable
examples and extend some of the results for the use of actual time series
analysis.

\section{Exactly solvable models}\label{sec:4}
In order to illustrate the convergence properties of EDMD,
analytic maps with accessible point spectrum are needed.
Although the Perron-Frobenius operator and its adjoint are
compact on Hardy-Hilbert spaces,
computing their eigenvalues remains a challenging task.
The first nontrivial family of analytic maps with explicitly
computable spectrum has been identified in \cite{SlBaJu_Nonl13}.
The family comprises circle maps $\tau$ which analytically extend to a neighbourhood of the entire
unit disk (not just an annulus), that is, maps arising from
Blaschke products. For these maps, the entire spectrum of the Perron-Frobenius operator is
determined by fixed point properties of $\tau$ inside the unit disk \cite{BaJuSl_AIHP17}.

To be more explicit consider a Blaschke product of
degree two, given by
\begin{equation}\label{ca}
\tau(z)=\frac{z-\mu}{1-\bar{\mu} z}\frac{z-\rho}{1-\bar{\rho} z},
\quad |\mu|, |\rho| <1,
\end{equation}
with two complex-valued parameters $\mu=|\mu| \exp(i\alpha)$
and $\rho=|\rho| \exp(i\beta)$, where $\bar{\mu}$ and $\bar{\rho}$ denote the complex conjugates
of $\mu$ and $\rho$, respectively. This map
preserves the unit circle, where (considered in angular coordinates)
it induces
a two-branch interval map
\begin{align}\label{cb}
\varphi \mapsto 2 \varphi &+
2 \mbox{arctan}\left(\frac{|\mu| \sin(\varphi-\alpha)}{1-|\mu|
\cos(\varphi-\alpha)}\right) \nonumber \\
&  +
2 \mbox{arctan}\left(\frac{|\rho| \sin(\varphi-\beta)}{1-|\rho|
\cos(\varphi-\beta)}\right) \quad (\mbox{mod } 2 \pi) \, .
\end{align}

The map \eqref{cb} can be considered as an analytic deformation of the
Bernoulli shift map which is obtained for the choice $\mu=\rho=0$.
The map $\tau$ is eventually expanding\footnote{
A map $\tau\colon\mathbb{T}\to\mathbb{T}$ is called eventually expanding if it has an iterate that is
expanding.} on the unit circle if and only if
\cite[Propositions 2.1, 3.1]{Pujals} it has a unique (attracting) fixed point $z_*=\tau(z_*)$ in the unit
disk, that is, $|z_*|<1$.
As shown\footnote{The results are stated for expanding Blaschke products, but can be extended to
the eventually expanding case.} in \cite{BaJuSl_AIHP17}, the
powers of the multiplier $\tau'(z_*)$ and their complex conjugates
are precisely the eigenvalues $\lambda_n$ of the Perron-Frobenius operator
\begin{equation}\label{cd}
\lambda_0=1, \quad \lambda_{2n-1}=(\tau'(z_*))^n,  \quad \lambda_{2n}=
(\overline{\tau'(z_*)})^n, \quad n \in \mathbb{N} \, .
\end{equation}

We use the map (\ref{ca}) to illustrate EDMD with a set of analytic
observables. An obvious choice is the set of the first $\bar{N}$
Fourier modes, that is, $\{\psi_k(z)=z^k : -\bar{N}\leq k \leq \bar{N}\}$ using
complex notation $z=\exp(i \varphi)$. Furthermore, as mentioned
in the previous section we evaluate \eqref{aba} and \eqref{abb}
for equidistant nodes on the unit circle
$z^{(m)}=\exp(2 \pi i m/M)$ with $m=0,\ldots,M-1$. It is straightforward to show that our
observables are orthogonal in the sense that $H_{k \ell}=\delta_{k,-\ell}$ in \eqref{abb},
if the number of nodes exceeds the number of observables,
$M\geq N=2\bar{N}+1$. It remains to evaluate \eqref{aba}, for which we will consider
maps of the form \eqref{ca}.

Let us first comment on the trivial parameter choice $\mu=\rho=0$, that is, on the
Bernoulli shift map. Clearly $\tau(z)=z^2$ has fixed point
$z_*=0$ with multiplier $\tau'(z_*)=0$, so all eigenvalues of the Perron-Frobenius operator
in \eqref{cd} apart from $\lambda_0$ vanish. For the application of EDMD, equation \eqref{aba}
can be easily evaluated to yield $G_{k \ell}=\delta_{-2k, \ell}$, as long as the number of
nodes is sufficiently large, that is, $M \geq 3 N/2$. Then the matrix
$A=G H^{-1}$ is given by $A_{k \ell}=\delta_{2k, \ell}$
and its eigenvalues in fact coincide with the leading part
of the exact spectrum given by \eqref{cd}.

\begin{figure}[!h]
\centering
\centering
\includegraphics[width=0.95\textwidth]{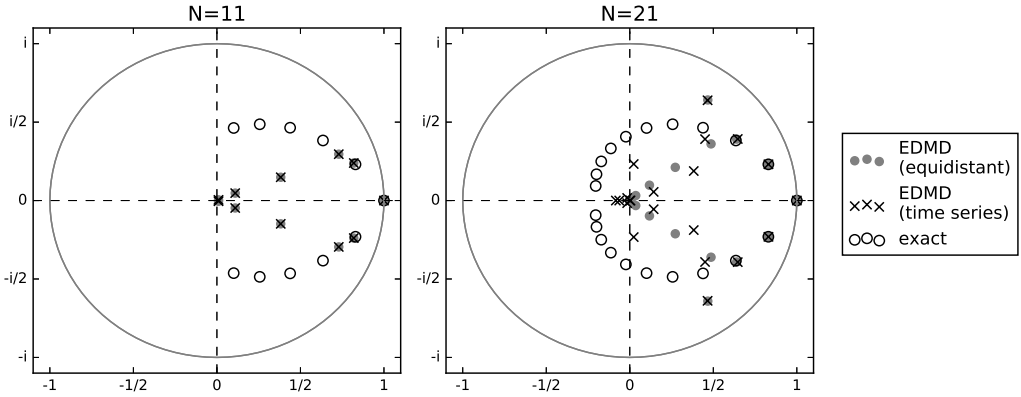}
\caption{Complex plane with exact eigenvalues (open symbols) and approximation by EDMD \corrn{with $M=100$ equidistant nodes} (full symbols)
\corrn{or by EDMD applied to time series of length $M=5\cdot 10^4$ (crosses)}
for the Blaschke product in \eqref{ca}, with \corrn{$\mu=\rho=0.33 \cdot \exp(i \pi/25)$} for
$N=11$ modes (left) or $N=21$ modes (right).}\label{fig1}
\end{figure}

For any non-trivial Blaschke map, the sums in \eqref{aba} and the related
finite-dimensional eigenvalue problem need to be evaluated numerically.
For that purpose we set \corrn{$\mu=\rho=0.33 \cdot \exp(i \pi/25)$}
which results in a spectrum with a fairly rich structure, so it
can serve as a test for the efficiency of EDMD.
Figure \ref{fig1} shows the eigenvalues of $A$ for $M=100$ equidistant nodes, $N=11$ and
$N=21$ observables compared with the exact expression \eqref{cd}. A set of $N=11$ modes
are just sufficient to approximate the subleading complex eigenvalue pair
($\lambda_1$ and $\lambda_2=\bar{\lambda}_1$) while all the other
values are spurious results. For a higher number of modes, $N=21$, about a quarter
of the eigenvalues of $A$ give reasonable estimates of the correct spectrum.
In particular, EDMD reproduces the leading part of the spectrum
of the compact Perron-Frobenius operator as asserted in Section \ref{sec:3}.

\begin{figure}[!h]
\centering
\includegraphics[width=0.7\textwidth]{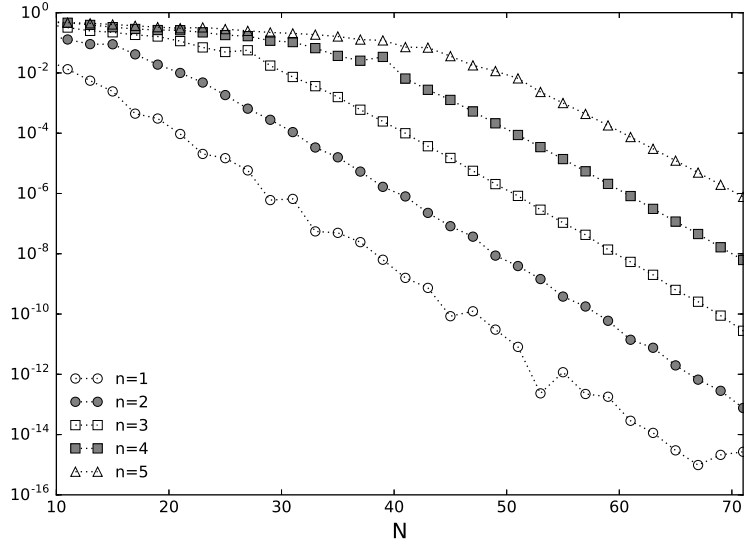}
\caption{Absolute error of the first five subleading complex-conjugate
eigenvalue pairs
$(\lambda_{2n-1}, \lambda_{2n})$ computed by EDMD
(with $M=1000$ nodes) for a Blaschke
product \eqref{ca}, with $\mu=\rho=0.33\cdot \exp(i\pi /25)$, as a function of the number of observables $N$.}\label{fig2}
\end{figure}

The dependence of the numerical error on the order $N$ of the \corr{eigenvalue} approximation
is shown in Figure \ref{fig2}. As we want to disentangle the effect of the two parameters, $M$ and $N$,
we take $M$ large enough, $M=1000$, so that all matrix elements
\corr{of the finite rank approximation $P_N\mathcal{L}P_N$} are estimated by the sums
\eqref{aba} sufficiently accurately. Hence any visible error is due to the finite mode
approximation. An exponential decay \corr{(in $N$)} of the \corr{eigenvalue approximation} error is clearly observable.

\section{Discussion}\label{sec:5}

We have shown that EDMD singles out eigenvalues of compact
Perron-Frobenius or Koopman operators arising in an analytic setting.
Moreover, the algorithm converges at an exponential rate in the number of observables used.\\

\textit{EDMD as a data analysis tool.}
So far we have focussed on equidistant points in phase space
to evaluate EDMD.
In applications one normally resorts to an actual time
series. If we reconsider the setup used for Figure \ref{fig1} but now take
nodes generated from a time series, that is, $z^{(m)}=z_m=\tau(z_{m-1})$,
we still obtain an accurate approximation of the spectrum of the Perron-Frobenius operator (see
crosses in Figure \ref{fig1}), as long as we minimise statistical fluctuations in the sums
\eqref{aba} and \eqref{abb} by taking a time series of sufficient length.
In fact, a slight modification of the arguments presented in Section
\ref{sec:3} allows one to base the matrix elements \eqref{bb} on integrals with respect to the analytic invariant density instead of the Lebesgue measure.  Whereas the matrix entries in \eqref{bb} change, the convergence
results remain unaffected\footnote{
  For a suitable density $\rho$,
  the appropriate
  Perron-Frobenius operator (with identical spectrum) is given by $\hat{\mathcal{L}}f = \rho^{-1}\mathcal{L}(f\rho)$, defined on a Hardy-Hilbert space with adapted inner product. The results in Section \ref{sec:3} then hold for
  $\hat{\mathcal{L}}$, with appropriately chosen basis and
  projection operators. In EDMD this change of basis is
  accounted for by the matrix $H^{-1}$ in \eqref{eq:A}.}.
  In particular, EDMD exhibits the same convergence, when used as a time series analysis tool.\\

\textit{Convergence in $M$ and $N$.}
Our approach is based on finite rank operators
represented by \eqref{bb}
converging in operator norm to a compact transfer operator. The nodes used in EDMD, in \eqref{aba} and \eqref{abb}, can naturally be considered as a sampling of the corresponding
integral.
However, a rigorous estimate which involves both
quantities, $M$ and $N$, appears to require taking the limit of large $M$ first. This
may in fact not be necessary.
While the approach so far was based on using orthogonal projectors,
one can in fact directly link the matrix representation involving sums with a compact transfer
operator by employing non-orthogonal projectors arising from collocation methods.
It may thus be possible to show
convergence for the case of $M$ and $N$ of the same order.\\

\begin{figure}[!h]
\centering
\includegraphics[width=0.4 \textwidth]{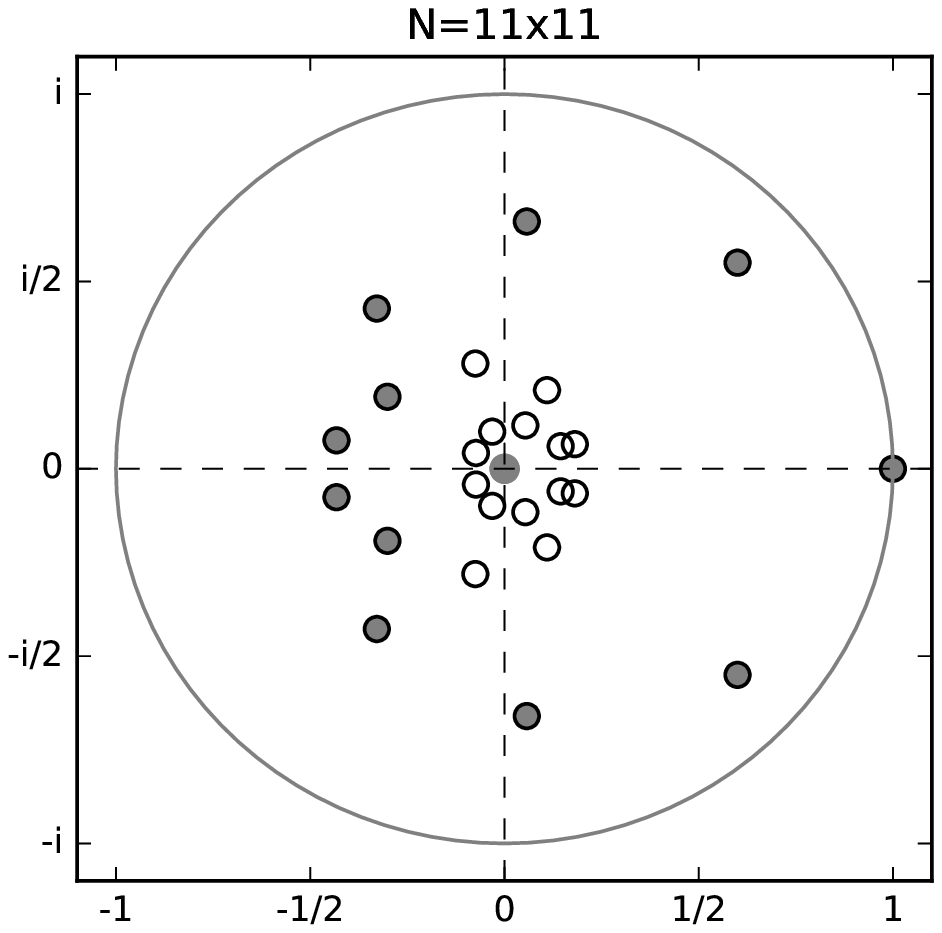}
\hspace{1.2cm}
\includegraphics[width=0.47 \textwidth]{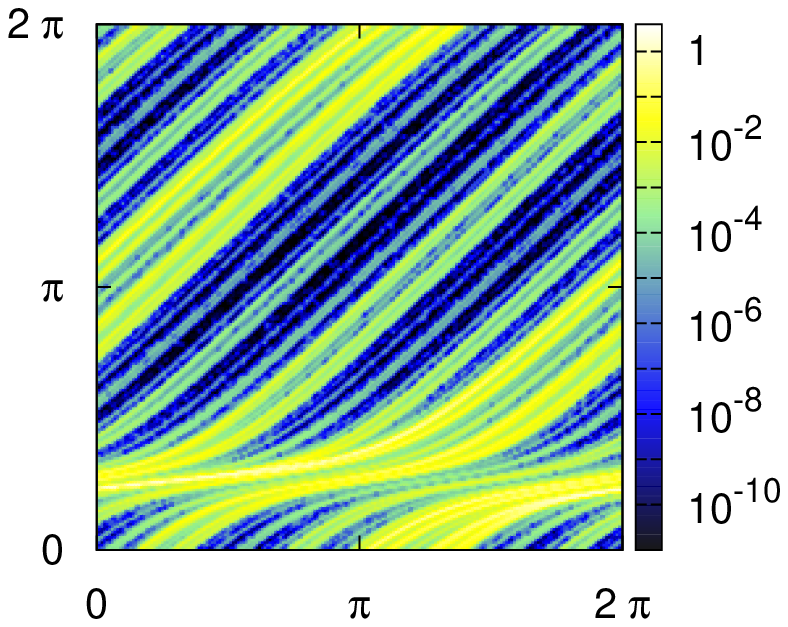}
\caption{Complex plane (left plot) with exact eigenvalues (open symbols) and approximation by
EDMD (full symbols) applied to an analytic deformation of the cat map \eqref{da} with \corrn{$
\mu=-0.6-0.55i$}.
Data obtained from EDMD with \corrn{$N=11 \times 11$} Fourier modes and
$M=201 \times 201$ nodes on a square lattice. Density plot (right) illustrating the invariant measure
for the map \eqref{da}
for the same parameter values.}\label{fig:fig4}
\end{figure}

\textit{Higher dimensional extensions and multifractal properties.}
So far one may object that we have enforced an analytic setting and that the results are
not really surprising as the maps do not allow for any complex multifractal behaviour. This
is in fact not correct as analyticity is only required for the actual
equation of motion, whereas the relevant invariant measure itself could be singular with
respect to the phase space volume. In order to
demonstrate this phenomenon, we resort to analytical solutions of two-dimensional
hyperbolic diffeomorphisms which
allow for fractal invariant measures if the Jacobian is not constant.
The presence of contracting and expanding directions requires using
more involved function spaces, that is, a particular class of anisotropic Hilbert spaces, for which
rigorous statements on spectral data of evolution operators are possible (see, for example,
\cite{SlBaJu_NONL17} for technical details).
We illustrate our point by considering an analytic deformation
of the cat map given by $(\varphi_1,\varphi_2)\mapsto (\varphi_1', \varphi_2')$ with
\begin{equation} \label{da}
\begin{aligned}
\varphi_1' =& 2 \varphi_1+\varphi_2 +
2 \mbox{arctan}\left(\frac{|\mu|
\sin(\varphi_1 + \varphi_2 -\alpha)}{1-|\mu|\cos(\varphi_1 + \varphi_2 -\alpha)}\right)\,,\\
\varphi_2' =& \hphantom{2}\varphi_1+\varphi_2 +
2 \mbox{arctan}\left(\frac{|\mu| \sin(\varphi_1 + \varphi_2 -\alpha)}{1-|\mu|
\cos(\varphi_1 + \varphi_2 -\alpha)}\right)\, ,
\end{aligned}
\end{equation}
where $\mu$ is a complex parameter with $|\mu|<1$.
This map is an analytic hyperbolic diffeomorphism of the torus, which for $\mu=0$ reduces to the
cat map $(\varphi_1,\varphi_2)\mapsto (2 \varphi_1+\varphi_2,\varphi_1+\varphi_2)$.
For non-vanishing $\mu$, the physical invariant measure is singular with respect to
phase volume with the corresponding invariant density exhibiting
fractal properties, see Figure \ref{fig:fig4} (right).
Employing more elaborate machinery one can show that
the corresponding Koopman operator is
compact on a suitable anisotropic Hilbert space (see \cite{SlBaJu_NONL17}
for similar results). Moreover, the
eigenvalues are determined by quantities associated with the fixed points of the map
\eqref{da}
in complex polydisks, and are of the form
$\lambda_1 = 0,
\lambda_{2n-1} = (-\mu)^{n}, \lambda_{2n} = \overline{(-\mu)}^{n},
 n\in\mathbb{N}$.
 Applying EDMD in this setting
 reproduces the leading part of the exact spectrum, see Figure \ref{fig:fig4} (left). Rigorous
 proofs of these statements will be presented elsewhere.

Altogether we have compelling evidence that convergence of EDMD points
towards an underlying compact operator structure which determines the
effective degrees of freedom in a complex dynamical setting.

\section*{Acknowledgement}
We gratefully acknowledge the support of the research presented in this article by the EPSRC
grant EP/R012008/1.

\end{document}